


\font\sc=cmcsc10 \rm


\newcount\secnb
\newcount\subnb
\newcount\parnb
\newcount\itemnb
\secnb=0
\subnb=0
\parnb=0
\itemnb=96
\newtoks\secref
\secref={}
\newtoks\subref
\subref={}

\def\smallskip{\par\vskip 2mm}
\def\medskip{\par\vskip 5mm}
\def\goodbreak{\penalty -100}


\def\section#1{\global\advance\secnb by 1
\secref=\expandafter{\the\secnb.}
\subref={}
\subnb=0\parnb=0\itemnb=96
\medskip\goodbreak
\centerline{\bf\S\the\secref\ #1}
\smallskip\nobreak}

\def\subsection#1{\global\advance\subnb by 1
\secref=\expandafter{\the\secnb.}
\subref=\expandafter{\the\subnb.}
\parnb=0\itemnb=96
\smallskip\goodbreak
\leftline{\bf\the\secref\the\subref\ #1}
\par\nobreak}

\def\references{\medskip\goodbreak
\centerline{\bf R\'ef\'erences}
\smallskip\nobreak}

\def\tit{\global\advance\parnb by 1
\itemnb=96
\smallskip\goodbreak
\noindent\the\secref\the\subref\the\parnb)\ }




\def\label#1{\relax}

\def\ref#1{\csname crossref#1\endcsname}


\newcount\bibrefnb
\bibrefnb=0

\def\biblabel#1{\global\advance\bibrefnb by 1\item{\the\bibrefnb.}}

\def\refto#1{\setbox1=\hbox{\sc #1}\copy1}
\def\reftosame{\hbox to \the\wd1{\hrulefill}}

\def\bibref#1{\csname bibref#1\endcsname}


\def\bye{\medskip\noindent{\it Addresse:} Benoit Fresse,\par
Laboratoire J.A. Dieudonn\'e,\par
Universit\'e de Nice,\par
Parc Valrose,\par
F-06108 Nice Cedex 02 (France).\par
\noindent{\it Courriel:} fresse@math.unice.fr\par
\vfill\end}

\def\beginresume{\smallskip
\begingroup\leftskip=1cm\rightskip=1cm
{\sc R\'esum\'e}. --- }
\def\endresume{\endgroup}
\def\beginabstract{\smallskip
\begingroup\leftskip=1cm\rightskip=1cm
{\sc Abstract}. --- }
\def\endabstract{\endgroup}


\def\Hom{\mathop{\rm Hom}\nolimits}
\def\O{{\cal O}}
\def\P{{\cal P}}
\def\C{{\bf C}}
\def\Z{{\bf Z}}

\def\Spec{\mathop{\rm Spec}\nolimits}
\def\Pois{{\mathop{{\cal P}ois}\nolimits}}
\def\can{{\mathop{can}\nolimits}}
\def\HH{\mathop{HH}\nolimits}


\hsize=15cm
\vsize=20.5cm

\headline={\ifnum\pageno=1{\it Projet de note aux C.R.Acad.Sci.Paris S\'er. I Math. (Alg\`ebre)}\hfil\else\hfil\fi}

\vbox{\medskip
\centerline{\bf Structures de Poisson sur une intersection compl\`ete \`a singularit\'es isol\'ees}
\smallskip
\centerline{\sc Benoit Fresse}
\smallskip
\centerline{14/2/2002}
\medskip
\beginresume
On \'etudie les structures de Poisson sur des vari\'et\'es singuli\`eres.
On consid\`ere dans ce but le complexe de Koszul associ\'e aux \'equations d'une intersection compl\`ete.
Ce complexe forme une alg\`ebre diff\'erentielle gradu\'ee qui est \'equivalente \`a l'alg\`ebre de la vari\'et\'e.
On montre qu'une structure de Poisson est \'equivalente \`a la donn\'ee d'une famille de multid\'erivations
sur le complexe de Koszul.
Si notre vari\'et\'e a des singularit\'es isol\'ees,
alors peut construire une telle famille de multid\'erivations de forme r\'eduite.
\endresume
\medskip
\centerline{\bf Poisson structures over a complete intersection with isolated singularities}
\smallskip
\beginabstract
We study Poisson structures over singular varieties.
In this purpose, we consider the Koszul complex associated to the equations of a complete intersection.
This complex forms a differential graded algebra which is equivalent to the algebra of the variety.
We show that a Poisson structure is equivalent to a sequence of multiderivations over the Koszul complex.
If our variety has isolated singularities,
then we can construct a sequence of multiderivations of reduced form.
\endabstract
\smallskip}

\medskip
\noindent{\bf Abridged English Version}
\smallskip

Suppose given a commutative algebra $\O$.
This algebra is equipped with a Poisson structure
if we have a biderivation $\{-,-\}: \O\otimes\O\,\rightarrow\,\O$
such that $\{f,\{g,h\}\} + \{g,\{h,f\}\} + \{h,\{f,g\}\} = 0$, for all $f,g,h\in\O$.
We would like to study such structures in the case $\O$ is the algebra
of a complete intersection with isolated singularities.
Our idea is to replace the algebra of a singular affine scheme
by an equivalent free differential graded algebra.

\tit{\it Differential graded modules}. ---
We work in the category of differential graded $\C$-modules (dg-modules, for short).
We adopt the following classical conventions.
A dg-module $V$ is either lower graded $V = V_*$ or upper graded $V = V^*$.
The relation $V^d = V_{-d}$ makes a lower grading equivalent to an upper grading.
The notation $|v|$ refers to the lower degree of a homogeneous element $v\in V$.
The differential of $V$ is denoted by $\delta: V\,\rightarrow\,V$.
The homology of the complex associated to $V$ is denoted by $H_*(V)$.
A morphism of dg-modules $V\,\rightarrow\,W$ is quasi-iso if the induced morphism $H_*(V)\,\rightarrow\,H_*(W)$ is iso.

The tensor product of dg-modules $V$ and $W$ is the dg-module such that
$(V\otimes W)_n = \bigoplus_{*\in\Z} V_*\otimes W_{n-*}$.
The differential of a tensor $v\otimes w\in V\otimes W$
is given by the classical formula
$\delta(v\otimes w) = \delta(v)\otimes w + \pm v\otimes\delta(w)$,
where $\pm = (-1)^{|v|}$.
We have a symmetry isomorphism $V\otimes W\,\rightarrow\,W\otimes V$
defined by $c(v\otimes w) = \pm w\otimes v$,
where $\pm = (-1)^{|v|\cdot |w|}$.
In the sequel, the notation $\pm$ refers to the sign produced by a permutation of homogeneous elements,
whose explicit determination follows from the definition of the symmetry isomorphism.

The module of homogeneous morphisms from $V$ to $W$ is the upper graded dg-module
such that $\Hom^n(V,W) = \prod_{*\in\Z} \Hom(V_{*+n},W_{*})$.
The differential of a homogeneous morphism $f\in\Hom^*(V,W)$
is provided by the commutator of $f$ with the differentials of $V$ and $W$.
We have explicitly $\delta(f) = \delta f - \pm f\delta$, where $\pm = (-1)^{|f|}$.
A $0$-cocycle $f\in\Hom^0(V,W)$ is equivalent to a morphism of dg-modules $f: V\,\rightarrow\,W$.

\tit{\it Differential graded algebras}. ---
A commutative dg-algebra is a dg-module $\tilde{\O}$
equipped with a product $\tilde{\O}\otimes\tilde{\O}\,\rightarrow\,\tilde{\O}$
which is associative and commutative.
The commutativity relation involves a sign because of the definition of the symmetry isomorphism
$\tilde{\O}\otimes\tilde{\O}\,\rightarrow\,\tilde{\O}\otimes\tilde{\O}$.

By convention,
a polynomial algebra $\tilde{\O} = \C[\tilde{u}_1,\ldots,\tilde{u}_r]$
denotes a commutative dg-algebra
which, as a graded commutative algebra, is freely generated
by homogeneous elements $\tilde{u}_1,\ldots,\tilde{u}_r\in\tilde{\O}$.
Let us mention that the differential of $\tilde{\O}$ is determined by the differential of the generators
$\delta(\tilde{u}_1),\ldots,\delta(\tilde{u}_r)\in\tilde{\O}$.

\tit{\it Resolutions}\label{EModel}. ---
A resolution of a commutative algebra $\O$ is a polynomial dg-algebra $\tilde{\O}$
such that $H_0(\tilde{\O}) = \O$ and $H_*(\tilde{\O}) = 0$ if $*>0$.
Equivalently, we assume that $\O$ is a dg-algebra concentrated in degree $0$.
The resolution $\tilde{\O}$ is a dg-algebra equipped with an augmentation morphism
$\epsilon: \tilde{\O}\,\rightarrow\,\O$
which is quasi-iso.

Let us assume that $\O$ is the algebra of a complete intersection
$h_1(x_1,\ldots,x_n) = \cdots = h_m(x_1,\ldots,x_n) = 0$.
In this case, a resolution of $\O$ is provided by the Koszul complex of the sequence $(h_1,\ldots,h_m)$,
which is nothing but the dg-algebra $\tilde{\O} = \C[x_1,\ldots,x_n,\tilde{h}_1,\ldots,\tilde{h}_m]$,
where $|\tilde{h}_j| = 1$ and $\delta(\tilde{h}_j) = h_j(x_1,\ldots,x_n)$.

\tit{\it On K\"ahler differentials}. ---
The definition of the algebra of K\"ahler differentials $\Omega^*(\tilde{\O})$ can be extended
to the case of a dg-algebra $\tilde{\O}$.
Precisely, we consider the dg-algebra over $\tilde{\O}$
generated by the elements $df$, where $f\in\tilde{\O}$,
together with the relations $d(f g) = df\cdot g + \pm f\cdot dg$,
where $\pm = (-1)^{|f|}$.
We have also $|df| = 1+|f|$ and $\delta(df) = - d(\delta(f))$.
(The differential symbol $d$ is equivalent to a tensor of degree $1$.)
The homogeneous component $\Omega^s(\tilde{\O})$ is the dg-module generated by the products
$f_0\cdot df_1\cdot\cdots\cdot df_s\in\Omega^*(\tilde{\O})$
(with $s$ differentials).

Similarly, we have a dg-module of multivectors $T^s(\tilde{\O})^*$
which is defined by the identity $T^s(\tilde{\O})^m = \Hom^n_{\tilde{\O}}(\Omega^s(\tilde{\O}),\tilde{\O})$.
An element $\tilde{\pi}\in T^s(\tilde{\O})^m$
is equivalent to an antisymmetric multiderivation
$\tilde{\pi}: \tilde{\O}\otimes\cdots\otimes\tilde{\O}\,\rightarrow\,\tilde{\O}$
which decreases the degree by $m-s$.
The superscript $s$ (respectively, $m$) refers to the order (respectively, to the degree)
of the multivector $\tilde{\pi}$.

\tit{\it Poisson structures}. ---
The dg-module of multivectors is endowed with the Schouten-Nijenhuis bracket
$[-,-]: T^s(\tilde{\O})^m\otimes T^t(\tilde{\O})^n\,\rightarrow\,T^{s+t-1}(\tilde{\O})^{m+n-1}$
as in the classical differential calculus.
Furthermore, a Poisson structure is equivalent to a bivector $\tilde{\pi}_2\in T^2(\tilde{\O})^2$
such that $\delta(\tilde{\pi}_2) = 0$ and $[\tilde{\pi}_2,\tilde{\pi}_2] = 0$.
We define a homotopy Poisson structure
as a sequence of multivectors $\tilde{\pi}_s\in T^s(\tilde{\O})^2$, $s\geq 2$,
whose sum $\tilde{\pi}_* = \sum_s\tilde{\pi}_s$ verifies the Maurer-Cartan equation
$\delta(\tilde{\pi}_*) + 1/2\cdot [\tilde{\pi}_*,\tilde{\pi}_*] = 0$.

\smallskip
We have obtained the following results:

\tit{\sc Theorem}\label{EHoPoisson}. --- {\it Let $\O$ be a commutative algebra over $\C$ equipped with a Poisson structure.
Let us fix a free differential graded resolution of $\O$ as in paragraph (\ref{EModel}).
The algebra $\tilde{\O}$ has a homotopy Poisson structure $\tilde{\pi}_*\in T^*(\tilde{\O})$
such that $\epsilon(\tilde{\pi}_2(df,dg)) = \{\epsilon(f)),\epsilon(g)\}$,
for all $f,g\in\tilde{\O}$.}

\tit{\sc Theorem}. --- {\it Let $\O$ be a commutative algebra over $\C$ equipped with a Poisson structure.
We assume that $\O$ is the algebra of a complete intersection
with isolated singularities $h_1(x_1,\ldots,x_n) = \cdots = h_m(x_1,\ldots,x_n) = 0$.
We consider the Koszul resolution of $\O$ introduced in the paragraph (\ref{EModel}).
We localize these algebras at $\P\in\Spec\O$.
There is a germ of homotopy Poisson structure $\tilde{\pi}_*\in T^*(\tilde{\O}_\P)$,
as in the theorem (\ref{EHoPoisson}),
such that, for $2\leq s\leq n-m$, we have
$$\tilde{\pi}_s(df_1,\ldots,df_s) = \left\{\matrix{ \tilde{\varpi}_s(df_1,\ldots,df_s,dh_1,\ldots,dh_m),
\qquad\hbox{if}\ \deg(f_1) = \cdots = \deg(f_s) = 0,\hfill\cr
0,\qquad\hbox{otherwise},\hfill\cr }\right.$$
where $\tilde{\varpi}_s: \Omega^{s+m}(\C[x_1,\ldots,x_n]_\P)\,\rightarrow\,(\tilde{\O}_\P)_{s-2}$.}

\medskip
\hbox to 3cm{\hrulefill}
\parnb=0\itemnb=0

\medskip
Soit $\O$ une alg\`ebre commutative sur $\C$.
Une structure de Poisson sur $\O$ est la donn\'ee d'une bid\'erivation $\{-,-\}: \O\otimes\O\,\rightarrow\,\O$
v\'erifiant la relation
$\{f,\{g,h\}\} + \{g,\{h,f\}\} + \{h,\{f,g\}\} = 0$,
quelque soient $f,g,h\in\O$.
On voudrait \'etudier ces structures quand $\O$ est l'alg\`ebre d'une intersection compl\`ete \`a singularit\'es isol\'ees.
L'id\'ee consiste \`a remplacer l'alg\`ebre d'un sch\'ema affine singulier
par une alg\`ebre diff\'erentielle gradu\'ee libre \'equivalente.

\section{Calcul diff\'erentiel gradu\'e}

\tit{\it Modules diff\'erentiels gradu\'es}. ---
On travaille dans la cat\'egorie des modules diff\'erentiels gradu\'es sur $\C$
(pour abr\'eger, on parlera de dg-modules).
On adopte les conventions classiques suivantes.
Un dg-module $V$ est soit gradu\'e inf\'erieurement $V = V_*$
soit gradu\'e sup\'erieurement $V = V^*$.
La relation $V^d = V_{-d}$ rend une graduation inf\'erieure \'equivalente \`a une relation sup\'erieure.
La notation $|v|$ renvoie au degr\'e d'un \'el\'ement homog\`ene $v\in V$.
La diff\'erentielle de $V$ est not\'ee $\delta: V\,\rightarrow\,V$.
L'homologie du complexe associ\'e \`a $V$ est not\'ee $H_*(V)$.
Un morphisme de dg-modules $V\,\rightarrow\,W$ est quasi-iso
si le morphisme induit $H_*(V)\,\rightarrow\,H_*(W)$ est iso.

Le produit tensoriel des dg-modules $V$ et $W$ est le dg-module
tel que $(V\otimes W)_n = \bigoplus_{*\in\Z} V_*\otimes W_{n-*}$.
La diff\'erentielle d'un tenseur $v\otimes w\in V\otimes W$ est donn\'ee par la formule classique
$\delta(v\otimes w) = \delta(v)\otimes w + \pm v\otimes\delta(w)$,
o\`u $\pm = (-1)^{|v|}$.
On a un isomorphisme de sym\'etrie $V\otimes W\,\rightarrow\,W\otimes V$
d\'efini par $c(v\otimes w) = \pm w\otimes v$,
o\`u $\pm = (-1)^{|v|\cdot |w|}$.
En g\'en\'eral,
la notation $\pm$ renvoie au signe produit par une permutation d'\'el\'ements homog\`enes,
dont la valeur r\'esulte de la d\'efinition de l'isomorphisme de sym\'etrie.

Le module des morphismes homog\`enes de $V$ dans $W$ est le dg-module gradu\'e sup\'erieurement
tel que $\Hom^n(V,W) = \prod_{*\in\Z}\Hom(V_{*+n},W_*)$.
La diff\'erentielle d'un morphisme homog\`ene $f\in\Hom^*(V,W)$ est donn\'ee par le commutateur de $f$
avec les diff\'erentielles de $V$ et $W$.
On a explicitement $\delta(f) = \delta f - \pm f \delta$,
avec $\pm = (-1)^{|f|}$.
Un $0$-cocycle $f\in\Hom^0(V,W)$ est \'equivalent \`a un morphisme de dg-modules $f: V\,\rightarrow\,W$.

\tit{\it Alg\`ebres diff\'erentielles gradu\'ees}. ---
Une dg-alg\`ebre commutative est un dg-module $\tilde{\O}$
muni d'un produit $\tilde{\O}\otimes\tilde{\O}\,\rightarrow\,\tilde{\O}$
qui est associatif est commutatif.
La relation de commutativit\'e comporte un signe d\^u \`a la d\'efinition de l'isomorphisme de sym\'etrie
$\tilde{\O}\otimes\tilde{\O}\,\rightarrow\,\tilde{\O}\otimes\tilde{\O}$.

Par convention, une alg\`ebre polynomiale $\tilde{\O} = \C[\tilde{u}_1,\ldots,\tilde{u}_r]$
d\'esigne une dg-alg\`ebre commutative qui, en tant que dg-alg\`ebre, est librement engendr\'ee
par des \'el\'ements homog\`enes $\tilde{u}_1,\ldots,\tilde{u}_r\in\tilde{\O}$.
La diff\'erentielle de $\tilde{\O}$ est d\'etermin\'ee par la diff\'erentielle des g\'en\'erateurs
$\delta(\tilde{u}_1),\ldots,\delta(\tilde{u}_r)\in\tilde{\O}$.

\tit{\it R\'esolutions}\label{Model}. ---
Une r\'esolution d'une alg\`ebre commutative $\O$ est une alg\`ebre polynomiale $\tilde{\O}$
telle que $H_0(\tilde{\O}) = \O$ et $H_*(\tilde{\O}) = 0$ si $* > 0$.
De fa\c con \'equivalente, on suppose que $\O$ est une dg-alg\`ebre concentr\'ee en degr\'e $0$.
La r\'esolution $\tilde{\O}$ est une dg-alg\`ebre munie d'un morphisme d'augmentation
$\epsilon: \tilde{\O}\,\rightarrow\,\O$ qui est quasi-iso.

Supposons par exemple que $\O$ est l'alg\`ebre d'une intersection compl\`ete
$h_1(x_1,\ldots,x_n) = \cdots = h_m(x_1,\ldots,x_n) = 0$.
On a alors une r\'esolution de $\O$ qui est fournie par le complexe de Koszul de la suite $(h_1,\ldots,h_m)$.
Ce complexe n'est rien d'autre que l'alg\`ebre diff\'erentielle gradu\'ee
$\tilde{\O} = \C[x_1,\ldots,x_n,\tilde{h}_1,\ldots,\tilde{h}_m]$
o\`u $|\tilde{h}_j| = 1$ et $\delta(\tilde{h}_j) = h_j(x_1,\ldots,x_n)$.
On note simplement que $\tilde{h}_i \tilde{h}_j = - \tilde{h}_j \tilde{h}_i$,
d'apr\`es la r\`egle de commutation des variables homog\`enes.

\section{Structures de Poisson homotopiques}

\tit{\it Formes de K\"ahler}. ---
La d\'efinition de l'alg\`ebre des formes de K\"ahler s'\'etend au cadre diff\'erentiel gradu\'e.
Ainsi, on d\'efinit $\Omega^*(\tilde{\O})$ comme la dg-alg\`ebre sur $\tilde{\O}$
engendr\'ee par les \'el\'ements $df$, avec $f\in\tilde{\O}$,
et quotient\'ee par les relations $d(f g) = df\cdot g + \pm f\cdot dg$,
o\`u $\pm = (-1)^{|f|}$.
On a aussi $|df| = |f| + 1$ et $\delta(df) = - d(\delta(f))$.
(Le symbole de diff\'erentielle $d$ est \'equivalent \`a un tenseur de degr\'e $1$.)
La composante homog\`ene $\Omega^s(\tilde{\O})$
est le dg-module engendr\'e par les produits
$f_0\cdot df_1\cdot\cdots\cdot df_s\in\Omega^*(\tilde{\O})$
(comportant $s$ diff\'erentielles).

\tit{\it Multivecteurs}. ---
On a \'egalement un dg-module de multivecteurs $T^s(\tilde{\O})^*$ d\'efini par la relation
$T^s(\tilde{\O})^m = \Hom_{\tilde{\O}}^m(\Omega^s(\tilde{\O}),\tilde{\O})$.
Un \'el\'ement $\tilde{\pi}\in T^s(\tilde{\O})^m$ \'equivaut \`a une multid\'erivation antisym\'etrique
$\tilde{\pi}: \tilde{\O}\otimes\cdots\otimes\tilde{\O}\,\rightarrow\,\tilde{\O}$
qui diminue le degr\'e de $m-s$.
L'exposant $s$ (respectivement, $m$) renvoie \`a l'ordre (respectivement, au degr\'e) du multivecteur $\tilde{\pi}$.

On munit le dg-module des multivecteurs du crochet de Schouten-Nijenhuis
$$[-,-]: T^s(\tilde{\O})^m\otimes T^t(\tilde{\O})^n\,\rightarrow\,T^{s+t-1}(\tilde{\O})^{m+n-1}$$
comme dans le calcul diff\'erentiel classique.
On a explicitement $[P,Q] = P\circ Q + \pm Q\circ P$,
o\`u $P\circ Q(df_1,\ldots,df_{s+t-1}) = \sum\pm P(df_{i_1},\ldots,df_{i_{s-1}},dQ(df_{j_1},\ldots,df_{j_t}))$,
pour tout multivecteurs $P\in T^s(\tilde{\O})^*$ et $Q\in T^t(\tilde{\O})^*$.
Les signes sont d\'etermin\'es par la r\`egle de commutation des tenseurs.
La somme s'\'etend sur l'ensemble des d\'ecompositions
$\{1,\ldots,s+t-1\} = \{i_1,\ldots,i_{s-1}\}\cup\{j_1,\ldots,j_t\}$.

\tit{\it Structures de Poisson}. ---
Une structure de Poisson homotopique est une suite de multivecteurs
$\tilde{\pi}_s\in T^s(\tilde{\O})^2$, $s\geq 2$,
dont la somme $\tilde{\pi}_* = \sum_s \tilde{\pi}_s$ (d\'efinie formellement)
v\'erifie l'\'equation de Maurer-Cartan
$\delta(\tilde{\pi}_*) + 1/2\cdot [\tilde{\pi}_*,\tilde{\pi}_*] = 0$.
On obtient les \'equations que v\'erifient les multivecteurs $\tilde{\pi}_s$, $s\geq 2$,
en d\'eveloppant les composantes d'ordre homog\`ene de cette \'equation.
On a explicitement:
$$\delta(\tilde{\pi}_s) + {1\over 2}\cdot\sum_{p+q = s+1} [\tilde{\pi}_p,\tilde{\pi}_q] = 0,$$
pour tout $s\geq 2$.
En particulier, pour $s = 2$, on obtient l'\'equation $\delta(\tilde{\pi}_2) = 0$.
Pour $s = 3$, on obtient $\delta(\tilde{\pi}_3) + 1/2\cdot [\tilde{\pi}_2,\tilde{\pi}_2]= 0$.

Une structure de Poisson \'equivaut \`a une structure de Poisson homotopique telle que $\tilde{\pi}_s = 0$ pour $s\not=2$.
Consid\'erons les multid\'erivations
$\{-,\cdots,-\}_s: \tilde{\O}\otimes\cdots\otimes\tilde{\O}\,\rightarrow\,\tilde{\O}$
associ\'ees aux multivecteurs $\tilde{\pi}_s\in T^s(\tilde{\O})^2$, $s\geq 2$.
On a par d\'efinition $\{f_1,\ldots,f_s\}_s = \pm\tilde{\pi}_s(df_1,\ldots,df_s)$.
Le signe provient de la commutation des symboles diff\'erentiels $d$
avec les \'el\'ements $f_1,\ldots,f_s$.
La composante d'ordre $s = 2$ de l'\'equation de Maurer-Cartan \'equivaut \`a la relation
$\delta\{f_1,f_2\}_2 = \{\delta(f_1),f_2\}_2 + \pm\{f_1,\delta(f_2)\}_2$.
La composante d'ordre $s = 3$ se r\'e\'ecrit:
$$\{\{f_1,f_2\}_2,f_3\}_2 + \pm\{\{f_2,f_3\}_2,f_1\}_2 + \pm\{\{f_3,f_1\}_2,f_2\}_2
= \pm\delta(\tilde{\pi}_3)(df_1,df_2,df_3).$$
Ainsi, en g\'en\'eral, on suppose seulement que la relation de Jacobi est v\'erifi\'ee modulo un cobord
(qui est repr\'esent\'e par le trivecteur $\tilde{\pi}_3$).

On a en fait d\'efini un type particulier de structure de Poisson homotopique.
On renvoie le lecteur au travail de Ginot pour une notion de structure homotopique plus g\'en\'erale
dans le cadre analogue des alg\`ebres de Gerstenhaber (cf. [\bibref{G}]).
Ces structures interviennent dans les travaux de Kontsevich et Tamarkin
sur la formalit\'e des complexes de Hochschild (cf. [\bibref{KF}]).

\smallskip
On a le r\'esultat suivant:

\tit{\sc Th\'eor\`eme}\label{HoPoisson}. --- {\it Soit $\O$ une alg\`ebre commutative sur $\C$ munie d'une structure de Poisson.
On fixe une r\'esolution de $\O$ comme dans le paragraphe (\ref{Model}).
L'alg\`ebre $\tilde{\O}$ a une structure de Poisson homotopique $\tilde{\pi}_*\in T^*(\tilde{\O})$
telle que $\epsilon(\tilde{\pi}_2(df,dg)) = \{\epsilon(f),\epsilon(g)\}$,
pour tout $f,g\in\tilde{\O}$.}

\smallskip
On construit la suite des multivecteurs pas \`a pas.
On observe que le morphisme d'augmen\-ta\-tion $\epsilon: \tilde{\O}\,\rightarrow\,\O$
induit un quasi-isomorphisme
$$\Hom_{\tilde{\O}}^*(\Omega^s(\tilde{\O}),\tilde{\O})\,\rightarrow\,\Hom_{\tilde{\O}}^*(\Omega^s(\tilde{\O}),\O)$$
et que les modules
$\Hom_{\tilde{\O}}^*(\Omega^s(\tilde{\O}),\O)$ ont une composante de degr\'e $*=2$ nulle d\`es que $s>2$.
Il s'ensuit que les obstructions \`a l'existence des multivecteurs $\tilde{\pi}_s$ sont nulles.

\section{Cas des intersections compl\`etes}

\tit{\sc Th\'eor\`eme}. --- {\it Soit $\O$ une alg\`ebre commutative sur $\C$ munie d'une structure de Poisson.
On suppose que $\O$ est l'alg\`ebre d'une intersection compl\`ete \`a singularit\'es isol\'ees
$h_1(x_1,\ldots,x_n) = \cdots = h_m(x_1,\ldots,x_n) = 0$.
On consid\`ere la r\'esolution de Koszul de $\O$ introduite dans le paragraphe (\ref{Model}).
On localise ces alg\`ebres en $\P\in\Spec\O$.
Il existe un germe de structure de Poisson homotopique $\tilde{\pi}_*\in T^*(\tilde{\O}_\P)$,
comme dans le th\'eor\`eme (\ref{HoPoisson}),
tel que, pour $2\leq s\leq n-m$, on ait
$$\tilde{\pi}_s(df_1,\ldots,df_s) = \left\{\matrix{ \tilde{\varpi}_s(df_1,\ldots,df_s,dh_1,\ldots,dh_m),
\qquad\hbox{si}\ \deg(f_1) = \cdots = \deg(f_s) = 0,\hfill\cr
0,\qquad\hbox{sinon},\hfill\cr }\right.$$
o\`u $\tilde{\varpi}_s: \Omega^{s+m}(\C[x_1,\ldots,x_n]_\P)\,\rightarrow\,(\tilde{\O}_\P)_{s-2}$.}

\smallskip
On consid\`ere la famille des complexes de Koszul g\'en\'eralis\'es associ\'es \`a la matrice Jacobienne
$((h_i)'_j)$ des polyn\^omes $(h_1,\ldots,h_m)$ (cf. [\bibref{N}, Appendice C]).
On utilise que ces complexes sont acycliques (le lieu singulier de notre vari\'et\'e \'etant de codimension $n-m$).
On forme le produit tensoriel de ces complexes de Koszul g\'en\'eralis\'es
avec le complexe de Koszul de $(h_1,\ldots,h_m)$.
On obtient ainsi un bicomplexe.
On construit les morphismes
$\tilde{\varpi}_s: \Omega^{s+m}(\C[x_1,\ldots,x_n]_\P)\,\rightarrow\,(\tilde{\O}_\P)_{s-2}$
comme le bord de certains cycles dans ce bicomplexe.

\tit{\it Exemple des hypersurfaces}. ---
On suppose $m = 1$ et $\O = \C[x_1,\ldots,x_n]/(h)$.
On a alors $\tilde{\O} = \C[x_1,\ldots,x_n,\tilde{h}]$, avec $|\tilde{h}| = 1$ et $\delta(\tilde{h}) = h$.
Le th\'eor\`eme ci dessus entraine que le crochet de Poisson de $\O$ est donn\'e par une formule de la forme
$$\{f,g\} = \sum_{i,j,k} p_{i j k} f'_i g'_j h'_k.$$
(On prend simplement $p_{i j k} = \tilde{\varpi}_2(dx_i,dx_j,dx_k)$.)
On a aussi une relation de la forme
$$\{\{f_1,f_2\},f_3\} + \{\{f_2,f_3\},f_1\} + \{\{f_3,f_1\},f_2\}
= h\cdot\sum_{i,j,k,l} q_{i j k l} (f_1)'_i (f_2)'_j (f_3)'_k h'_l.$$
(On consid\'ere les polyn\^omes $q_{i j k l}$
tels que $\tilde{\varpi}_3(dx_i,dx_j,dx_k,dx_l) = \tilde{h}\cdot q_{i j k l}$
dans $\tilde{\O}_1$.)

Comme application, on retrouve les formules que Alev et Lambre ont obtenues par un calcul direct
pour le crochet de Poisson symplectique des surfaces de Klein (cf. [\bibref{AL}]).
Les polyn\^omes $p_{i j k}$ sont alors n\'ecessairement constants pour des raisons d'homog\'en\'eit\'e.
On peut \'egalement supposer que les polyn\^omes $q_{i j k l}$ sont nuls dans ce cas.

\tit{\it Homologie de Poisson}. ---
Le th\'eor\`eme de structure ci-dessus nous permet d'obtenir des r\'esultats sur l'homologie de Poisson
des intersections compl\`etes \`a singularit\'es isol\'ees.
On consid\`ere l'homologie de Poisson $H^\Pois_*(\O,\O)$ telle qu'elle est d\'efinie
dans la th\'eorie des op\'erades de Koszul (cf. [\bibref{F}]).
L'homologie de Poisson d\'efinie par Koszul et Brylinski (cf. [\bibref{K}])
est not\'ee $H^\can_*(\O,\O)$
et sera d\'esign\'ee comme l'homologie canonique de $\O$.
On rappelle que l'on a une suite spectrale
$E^1_{s,t} = \HH^{(s)}_{s+t}(\O,\O)\,\Rightarrow\,H^\Pois_{s+t}(\O,\O)$,
o\`u $\HH^{(s)}_*$ d\'esigne la composante de poids $s$ de l'homologie de Hochschild.
On a obtenu le r\'esultat suivant:

\tit{\sc Th\'eor\`eme}. --- {\it Si $\O$ est l'alg\`ebre d'une intersection compl\`ete \`a singularit\'es isol\'ees,
alors la suite spectrale $E^1_{s,t} = HH^{(s)}_{s+t}(\O,\O)\,\Rightarrow\,H^\Pois_{s+t}(\O,\O)$
d\'eg\'en\`ere au rang $2$.
On a de plus $H^\Pois_*(\O,\O) = H^\can_*(\O,\O)$ en degr\'e $*\leq\dim\O$
et $H^\Pois_*(\O,\O) = \HH_*(\O,\O)$ en degr\'e $*\geq\dim\O$.}

\smallskip
On utilise que les composantes de poids de l'homologie de Hochschild
sont d\'etermin\'ees par l'homologie des complexes de Koszul g\'en\'eralis\'es
introduits pr\'ec\'edemment (cf. [\bibref{LR}]).

\smallskip
On peut \'enoncer un r\'esultat plus pr\'ecis dans le cas d'une hypersurface:

\tit{\sc Th\'eor\`eme}. --- {\it On suppose que $\O$ est l'alg\`ebre
d'une hypersurface \`a singularit\'es isol\'ees $h(x_1,\ldots,x_n) = 0$
telle que $h\in(h'_1,\ldots,h'_n)$.
On a alors $H^\Pois_*(\O,\O) = H^\can_*(\O,\O)$ pour $*\leq n-1$
et $H^\Pois_*(\O,\O) = \O/(h'_1,\ldots,h'_n)$ pour $*\geq n-1$.}

\references{\parindent=0cm\leftskip=1cm\rightskip=1cm

\biblabel{AL} \refto{J. Alev, T. Lambre},
{\it Comparaison de l'homologie de Hochschild et de l'homologie de Poisson
pour une d\'eformation des surfaces de Klein},
in ``Algebra and operator theory, Tashkent, 1997'',
Kluwer Acad. Publ. (1998), 25-38.

\biblabel{F} \refto{B. Fresse}, {\it Homologie de Quillen pour les alg\`ebres de Poisson},
C. R. Acad. Sci. Paris S\'er. I Math. {\bf 326} (1998), 1053-1058.

\biblabel{G} \refto{G. Ginot}, {\it Homologie et mod\`ele minimal des alg\`ebres de Gerstenhaber},
pr\'epublication (2001).

\biblabel{KF} \refto{M. Kontsevich}, {\it Operads and motives in deformation quantization},
Lett. Math. Phys. {\bf 48} (1999), 35-72. 

\biblabel{K} \refto{J.-L. Koszul},
{\it Crochet de Schouten-Nijenhuis et cohomologie},
in ``Elie Cartan et les math\'ematiques d'aujourd'hui'',
Ast\'erisque hors s\'erie (1985), 257-271.

\biblabel{LR} \refto{A. Lago, A.G. Rodicio},
{\it Generalized Koszul complexes and Hochschild (co)-ho\-mo\-lo\-gy of complete intersections},
Invent. Math. {\bf 107} (1992), 433-446.

\biblabel{N} \refto{D.G. Northcott}, Finite free resolutions,
Cambridge tracts in mathematics {\bf 71}, Cambridge university press, 1976.

}

\bye